\documentclass[a4paper,11pt]{article}
\usepackage{cancel}
\usepackage{amsmath}
\usepackage{amssymb}
\usepackage{array}
\usepackage{ragged2e}
\usepackage{algpseudocode}

\usepackage{soul}
\usepackage{rotating}
\usepackage{adjustbox}
\usepackage{bm}
\usepackage{longtable} 
\usepackage{arydshln}
\usepackage{adjustbox, lipsum}
\usepackage{enumerate}
\usepackage{rotating}
\usepackage[utf8]{inputenc}
\usepackage{amsfonts}
\usepackage{pstricks}
\usepackage{pdflscape}
\usepackage{graphicx}
\usepackage{soul}
\usepackage{afterpage}
\usepackage{graphics}
\usepackage{adjustbox}
\usepackage{rotating}
\usepackage{todonotes}

\usepackage[breaklinks]{hyperref}

\usepackage{cleveref}

\usepackage{tikz}
\usepackage{pgfplots}
\usepackage{rotating}
\usepackage[numbers,sort]{natbib}
\usepackage{booktabs}
\usepackage{setspace}

\usepackage{authblk}
\usepackage{multirow}
\usepackage{mathtools}
\usepackage{scalerel}
\usepackage{algorithm}%
\usepackage{algorithmicx}%

\newtheorem{prop}{\bf Proposition}[section]

\newtheorem{rem}{\bf Remark}[section]

\DeclarePairedDelimiter\abs{\lvert}{\rvert}%

\newcommand{\dem}{\par \noindent{\bf Proof:} }
\newcommand{\fin}{\hfill $\square$  \par \bigskip}

\usepackage{xcolor}
\definecolor{gr}{RGB}{0,153,0}

\usepackage{ulem}
\title{Constraint relaxation for the Discrete Ordered Median Problem}
\date{}

\author[a,b]{Luisa I. Martínez-Merino}
\author[b,c]{Diego Ponce \footnote{Corresponding author: dponce@us.es (D. Ponce)}}
\author[b,c]{Justo Puerto}

\affil[a]{\small {Departamento de Estad\'istica e Investigaci\'on Operativa, Universidad de C\'adiz}}
\affil[b]{\small {Instituto de Matemáticas de la Universidad de Sevilla (IMUS)}}
\affil[c]{\small{Departamento de Estadística e Investigación Operativa, Universidad de Sevilla}}

\allowdisplaybreaks
\setcitestyle{authoryear,open={(},close={)}}

\makeatletter
\let\oldabs\abs
\def\abs{\@ifstar{\oldabs}{\oldabs*}}
\makeatother

\begin{document}

\maketitle
	
\begin{abstract}This paper compares different exact approaches to solve the Discrete Ordered Median Problem (DOMP).
		In recent years, DOMP has been formulated using set packing constraints giving rise to one of its most promising formulations. The use of this family of constraints, known as \emph{strong order constraints} (SOC), has been validated in the literature by its theoretical properties and because their linear relaxation provides very good lower bounds. Furthermore, embedded in branch-and-cut or branch-price-and-cut procedures as valid inequalities, they allow one to improve computational aspects of solution methods such as CPU time and use of memory. In spite of that, the above mentioned  formulations require to include another family of  order constraints, e.g., the \emph{weak order constraints} (WOC), which leads to coefficient matrices with elements other than \{0,1\}. In this work, we develop a new approach that does not consider extra families of order constraints and furthermore relaxes SOC -in a branch-and-cut procedure that does not start with a complete formulation- to add them iteratively using row generation techniques to certify feasibility and optimality. Exhaustive computational experiments show that it is advisable to use row generation techniques in order to only  consider \{0,1\}-coefficient matrices modeling the DOMP. Moreover, we test how to exploit the problem structure. Implementing an efficient separation of SOC using callbacks improves the solution performance. This allows us to deal with bigger instances than using fixed cuts/constraints pools automatically added by the solver in the branch-and-cut for SOC, concerning both the formulation based on WOC and the row generation procedure.
\end{abstract}

	
\section{Introduction}

At times, very hard combinatorial optimization problems  contain easy combinatorial subproblems after relaxing some of their constraints. A paradigmatic example is the Traveling Salesman Problem: after the elimination of its subtour elimination constraints it turns into the Linear Assignment Problem which is polynomially solvable. This pattern calls for developing techniques that take advantage of this situation to solve some combinatorial problems based on their constraint relaxation. This approach is not new and the reader is referred to \cite{Focacci1999,Focacci2002a,Focacci2002b} and the references therein for further details.

This behavior is not only observed in problems where formulations include an exponential number of constraints. Actually, it also occurs in many polynomial size formulations. One of these cases is the Discrete Ordered Median Problem (DOMP) modeled with the strong order constraints formulation as introduced in \cite{Labbe2017}. If one removes the family of strong order constraints, whose acronym is SOC, the resulting problem is the standard $p$-median problem that is known to be combinatorially friendly  \citep{Hakimi1964,Marin2019,ReVelle1970}.
The aim of this work is to  develop  solution techniques for  DOMP based on  constraint relaxations.

The DOMP is  a discrete location model that allows to generalize several classical discrete location problems \citep[see, e.g.,][]{NickelPuerto2005}. For instance, the discrete $p$-center and $p$-median are particular cases of DOMP. 
Assume that we are given a set of clients, a set of candidate locations for facilities, and the allocation costs from each candidate facility to each client.	
The objective of DOMP is to locate $p$ facilities in such a way that a certain weighted function of the allocation costs is minimized. These weights are not assigned to specific costs but to their sorted values. Namely, the weighted average sorts the allocation costs in a nondecreasing order
and then, it performs the scalar product of this so obtained
sorted cost vector by the vector of weights.

In the literature, one can find different applications of the ordered median operator.  For instance, it has been applied to facility location \citep{Aouad2019,Dominguez2020,Espejo2009,Kalcsics2010, MM2017,Tamir2001}, multicast communication \citep{Fourour2020}, multiobjective Markov decision processes \citep{Ogryczak2011}, voting problems \citep{Ponce2018}, supervised classification \citep{Marin2022}, tomography reconstruction  \citep{Calvino2022}, and network design \citep{Puerto2013}, 
among other situations. 

DOMP was first introduced in \cite{Nickel2001} as an integer nonlinear problem. Then, in \cite{Boland2006}, this problem was modeled as a mixed integer linear program. 
Some works on this problem take advantage of some particular characteristics. Specifically, \cite{Marin2009} introduce an efficient covering formulation for DOMP considering free self-service, ties in the cost matrix and a non-negative weighted order vectors in the objective function. Futhermore, \cite{Marin2010} present a covering reformulation for weighted order vectors containing zeros and an extended model for vectors even with negative elements.

In \cite{Labbe2017}, a new three-index formulation based on set packing constraints is proposed. These set packing constraints are known as {\sl strong order constraints} or SOC, and the number of these constraints appearing in the formulation is $\mathcal{O}(n^3)$.
In addition, another new formulation, solved by an efficient branch-and-cut procedure that 
provides good results in terms of time, is introduced.

This second formulation is based on the aggregation of the SOC corresponding to the same position. The resulting order constraints are the {\sl weak order constraints} (from now on WOC). This formulation includes SOC as valid inequalities. 
Both formulations present small integrality gaps.

Recently, in \cite{Deleplanque2020}, a novel branch-price-and-cut algorithm has been proposed. This procedure is based on a formulation with an exponential number of variables that corresponds to a set partitioning model. The proposed approach 
allows to handle larger instances since it requires less memory to run the model.

In this paper, we want to explore different exact approaches to solve DOMP. The first one uses branch-and-cut techniques based on one of the most promising formulations, namely the formulation based on WOC proposed in \cite{Labbe2017}, adding SOC as valid inequalities. 
Additionally, we compare the use of cut pools in the branch-and-cut with respect to the use of callbacks to implement an \emph{ad hoc} separator proposed in \cite{Labbe2017}. By setting up pools of cuts, all SOC are initially stored and then, solvers decide which cuts are included during the branch-and-cut process. In contrast, applying the callbacks using the separation algorithm introduced in \cite{Labbe2017}, SOC are not initially stored and the implementation of the  SOC separation is based on a sequential update of the left-hand side of the corresponding order constraints. This separation can be performed in $\mathcal{O}(n^3)$. 

The second method is based on a constraint relaxation on the formulation using SOC to model the order.
This procedure, to solve DOMP, starts with a relaxed formulation where all SOC are removed and feasibility is enforced adding model constraints from the SOC family in the searching tree. Although the number of SOC is polynomial, $\mathcal{O}(n^3)$, this number of constraints becomes too large to be handled when $n$ increases. Consequently, it is interesting to study this approach since we could improve the time and memory performance by only including the necessary constraints in the solution process. Again, in this case, we compare the branch-and-cut procedure through callbacks with respect to the branch-and-cut based on constraint pools.

The contributions of this paper can be summarized as follows:
\begin{enumerate}
	\item Comparing a branch-and-cut approach to solve DOMP based on the so called WOC formulation with a constraint relaxation approach for DOMP based on removing SOC.
	\item Comparing the performance of the branch-and-cut and the constraint relaxation approach when using callbacks based on specific tailor made separation oracles with respect to the use of fixed pools of cuts/constraints.
	\item Reporting intensive computational tests which show the limits of the different considered solution methods.
\end{enumerate}

The remainder of this work is organized as follows.  In Section \ref{section:1}, we introduce the notation and description of DOMP. Besides, we recall two formulations for DOMP that will be used along the paper ($\text{DOMP}_{\text{WOC}}$ and $\text{DOMP}_{\text{SOC}}$). In Section \ref{section:3}, we present two solution methods for the problem. First, we describe the branch-and-cut procedure for $\text{DOMP}_{\text{WOC}}$  introduced in \cite{Labbe2017}. Then, we propose a novel row generation procedure for $\text{DOMP}_{\text{SOC}}$. Section \ref{computational} is devoted to the analysis of both solution methods. In addition, we present a comparison between the results of using pools of cuts/contraints and using callbacks for the branch-and-cut and the row generation techniques. Finally, in Section \ref{section:5}, we include some conclusions and future research lines.

\section{Problem definition and formulations}\label{section:1}

This section is devoted to recall the definition and some formulations of DOMP that will be instrumental in our discussion.	We shall follow the following notation. We denote by $I=\{1,\ldots,n\}$ the set of $n$ clients and, at the same time and without loss of generality, the set of $n$ potential facility locations. Facilities are assumed to be uncapacitated, i.e., they can supply as many clients as desired. Besides, $c_{ij}$ denotes the cost for serving client $i$ from facility $j$, for $i,j\in I$. 

Given a set $J$ composed by $p$ open facilities, $c_i(J)$ represents the cost of allocating client $i$ to the  facility set $J$, i.e., $ c_i(J)=\min_{j\in J}c_{ij}$. In addition, if the vector of costs $c_i(J)_{i=1,\dots,n}$ for $J\subset I$ is sorted in non-decreasing order, we denote by $c^{(k)}(J)$ the allocation cost in position $k\in K=I$ of this sorted vector. Thus, it holds that $c^{(1)}(J)\leq c^{(2)}(J)\leq\ldots\leq c^{(n)}(J)$.

The aim of DOMP is to determine a subset of $p$ facilities $J\subset I$  to open, and to assign each client to an open facility in order to minimize the ordered median objective function. Given a vector $\lambda=(\lambda^k)_{k\in K}$ such that $\lambda^k\geq 0$, $k\in K$, the objective function of DOMP can be expressed as 
$\sum_{k\in K}\lambda^k c^{(k)}(J).$ Consequently, the definition of DOMP is
\begin{equation}
\min_{J\subset I:\lvert J\rvert=p}\sum_{k\in K}\lambda^k c^{(k)}(J).\tag{DOMP}
\end{equation}

Observe that this formulation generalizes several standard discrete location problems. For instance: if $\lambda^1=\lambda^2=\ldots=\lambda^n=1$, this model is the $p$-median problem; if $\lambda^1=\lambda^2=\lambda^{n-1}=0,\lambda^n=1$, one gets the $p$-center; if $\lambda^1=\lambda^2=\ldots=\lambda^{n-k}=0,\lambda^{n-k+1}=\ldots=\lambda^n=1$, the resulting problem is the $k$-centrum; etc.

DOMP is known to be $\mathcal{NP}$-complete, see \cite{NickelPuerto2005}, and as mentioned in the introduction, different formulations have been proposed to deal with this problem. In this paper, we will elaborate on two of the most recent and promising formulations presented in \cite{Labbe2017}. On the one hand, that paper introduces a three-index formulation where order is modeled by a family of set packing constraints, SOC. On the other hand, it  also presents an aggregated version of that formulation where the order is ensured by a different set of constraints called WOC. We will build on those two formulations, thus in order to be self-contained, in the following subsection we provide full details of them.

\subsection{Strong Order Constraints formulation}
For the formulation based on {\sl strong order constraints}, the next families of variables are required:
\begin{eqnarray*}
	y_{j}&=&\begin{cases}
		1,&\mbox{if facility $j$ is open,}\\
		0,&\mbox{otherwise,}\\
	\end{cases} \mbox{ for $j\in I$,}\\
	x_{ij}^k&=&\begin{cases}
		1,&\mbox{if client $i$ is allocated to facility $j$ and the}\\
		&\mbox{associated cost is in position $k$ of the}\\
		& \mbox{sorted sequence of allocation costs,}\\
		0,&\mbox{otherwise},\\
	\end{cases}\mbox{for $i,j\in I$, $k\in K$.}
\end{eqnarray*}

Besides, we denote the rank of the allocation cost $c_{ij}$ by  $r_{ij}$, i.e.,  $r_{ij}=\ell$ if $c_{ij}$ is the $\ell$-th element in the  list of the costs  $c_{ij}$, for all $i,j \in I$, sorted in a  non decreasing sequence and where ties are broken arbitrarily. Then, $\text{DOMP}_{\text{SOC}}$ formulation is the following.

\begin{align}
\textbf{($\text{DOMP}_{\text{SOC}}$)}  \min&&\displaystyle\sum_{i\in I}\sum_{j\in I}\sum_{k\in K}\lambda^kc_{ij}x_{ij}^k\label{c3:ofdomp4}\\*
\mbox{s.t.}&&\displaystyle\sum_{j\in I}\sum_{k\in K}x_{ij}^k&=1,&&i \in I,\label{eq:1fb}\\*
&&\displaystyle\sum_{i\in I}\sum_{j\in I}x_{ij}^k &= 1,&&k\in K,\label{eq:2fb}\\*
&&\displaystyle\sum_{k\in K}x_{ij}^k &\leq y_j, && i,j\in I, \label{in:1fb}\\*
&&\displaystyle\sum_{j \in I}y_{j}& = p,\label{eq:3fb}\\*
&&\hspace*{-5cm}\sum_{i'\in I}\sum_{\substack{j'\in I:\\ r_{i'j'}\leq r_{ij}}}x_{i'j'}^k + \sum_{i'\in I}\sum_{\substack{j'\in I:\\ r_{i'j'}\geq r_{ij}}}x_{i'j'}^{k-1}  &\le 1, && i,j \in I,\; k\in K, k \neq 1,\label{in:3f1}\tag{\text{SOC}}\\*
&&x_{ij}^k,\; y_j   &\in \{0,1\},&& i,j \in I,k\in K. \label{binary:xf1}
\end{align}

Constraints \eqref{eq:1fb} ensure that each client is served by just one facility in one position. Similarly, constraints \eqref{eq:2fb} are necessary to guarantee that only one allocation cost is in each sorted position. Constraints \eqref{in:1fb} ensure that each client is allocated to an open facility and that the allocation cost of a client can only be in at most one position. Constraint \eqref{eq:3fb} restricts that exactly $p$ facilities must be open. Constraints \eqref{in:3f1} are the so-called {\sl strong order constraints} which ensure the correct sorting of the allocation costs. These constraints are set packing constraints, i.e., at most one of the variables of the l. h. s. could take value one. The incompatibility of two or more variables taking value one is due to \eqref{in:1fb} and the fact that a variable $x_{ij}^k$ cannot take value one if  $x_{i'j'}^{k-1}=1$ when $r_{ij}<r_{i'j'}$, for $i,j,i',j' \in I, k\in K, k \neq 1$. We refer the reader to \cite{Labbe2017} for a more detailed explanation. Finally, constraints \eqref{binary:xf1} are the domain of definition of the variables. The reader may note that removing \eqref{in:3f1} the formulation results in the $p$-median.

\subsection{Weak Order Constraints formulation}

Despite the good mathematical properties of $\text{DOMP}_{\text{SOC}}$, as the lower bound given by its linear relaxation, the number of \eqref{in:3f1} is $\mathcal{O}(n^3)$. Consequently, when the number of clients increases, the number of \eqref{in:3f1} becomes too large to be handled by a solver. For this reason, \cite{Labbe2017} introduce an alternative family of constraints to ensure the order of costs.  These new constraints are based on the aggregation of \eqref{in:3f1} corresponding to the same position. (The reader is referred to \cite{Labbe2017} for further details.) This alternative formulation results in the following.

\begin{align}
\textbf{($\text{DOMP}_{\text{WOC}}$)}  \min&&\displaystyle\sum_{i\in I}\sum_{j\in I}\sum_{k\in K}\lambda^kc_{ij}x_{ij}^k\nonumber\\*
\mbox{s.t.}&&\eqref{eq:1fb} - \eqref{binary:xf1},\nonumber\\*
&&\hspace*{-4cm}\sum_{i\in I}\sum_{j\in I}\left(\sum_{i'\in I}\sum_{\substack{j'\in I:\\ r_{i'j'}\leq r_{ij}}}x_{i'j'}^k + \sum_{i'\in I}\sum_{\substack{j'\in I:\\ r_{i'j'}\geq r_{ij}}}x_{i'j'}^{k-1}  \right)&\le n^2, && k\in K, k \neq 1.\label{in:3f2}\tag{\text{WOC}}
\end{align}

Constraints labeled by \eqref{in:3f2} are known as \emph{weak order constraints}. They ensure that if facility $j$ serves client $i$ and its cost $c_{ij}$ is in position $k$ of the sorted cost vector of the solution, then there must be a smaller or equal allocation cost in position $k-1$. This is due to the coefficients corresponding to each variable in the constraints. In each inequality, there are represented two positions ($k-1$ and $k$). By constraints \eqref{eq:2fb}, only two variables must take value 1, and the remaining ones take value 0. Assuming that the variables with value one for position $k$ and $k-1$ correspond to positions $s$ and $t$ of the sorted costs, respectively, the inequality can be expressed as follows: 
$$(n^2-(s-1))x_{i_s j_s}^k+tx_{i_t j_t}^{k-1}\leq n^2,$$
with $i_s,j_s,i_t,j_t\in I$ such that $c_{i_sj_s}$ and $c_{i_tj_t}$ are the $s$-th and $t$-th smallest allocation costs in matrix $(c_{ij})_{n\times n}$. This is valid if and only if $t<s$.

In \cite{Labbe2017}, it is shown that $\text{DOMP}_{\text{SOC}}$ formulation  provides a relevant improvement of the integrality gap with respect to $\text{DOMP}_{\text{WOC}}$ formulation. In other words, for most of the instances, fractional solutions satisfying \eqref{in:3f2} could be cut including \eqref{in:3f1}. Thus, it is recommended to use (SOC) as valid inequalities of $\text{DOMP}_{\text{WOC}}$.

\section{Solution methods}\label{section:3}

Both formulations,  $\text{DOMP}_{\text{SOC}}$ and $\text{DOMP}_{\text{WOC}}$, can be solved by using standard MIP solvers as CPLEX, Gurobi, Xpress, or SoPlex. However, the good performance of these formulations is rather limited for large sizes of the problem as we will see in Section \ref{computational}. The reader should observe that already for $n=100$ clients the number of (SOC) is almost $10^6$.

To improve the performance of $\text{DOMP}_{\text{WOC}}$, \cite{Labbe2017} propose a branch-and-cut procedure which starts by solving the linear program relaxation of $\text{DOMP}_{\text{WOC}}$, and then it includes \eqref{in:3f1} as valid inequalities whenever necessary. This procedure is described in Section \ref{section:B&C}. 

In this paper, we propose an alternative solution method which exploits the good properties of $\text{DOMP}_{\text{SOC}}$ avoiding the use of the complete family of strong order constraints. This solution method consists in a row generation procedure which initially considers $\text{DOMP}_{\text{SOC}}$ without \eqref{in:3f1}, and then it iteratively includes these order constraints. As far as we know, this row generation method has not been considered before for DOMP. Section \ref{section:rowg} is devoted to the description of this procedure.

\subsection{Branch-and-cut for $\text{DOMP}_{\text{WOC}}$}\label{section:B&C}

The branch-and-cut procedure has become an efficient method for solving large instances of models where the number of constraints is intractable for solvers. For instance, it has been successfully applied to the matching problem with blossoms \citep{Edmonds,Grtschel1985,Letchford}, problems related to trees \citep{FernandezPozoPuertoScozzari,Magnanti}, clustering \citep{Benati2017}, and the orienteering arc routing problem \citep{Archetti16}, to name a few. In \cite{Labbe2017}, a branch-and-cut method for \eqref{in:3f1} in  $\text{DOMP}_{\text{WOC}}$ formulation is proposed. 

This branch-and-cut procedure could be handled by two different perspectives. On the one hand, most of the current solvers have some options to define a fixed pool of cuts that are added automatically when necessary in the cut generation. This means that, changing some parameters in the solver,  we can remove some families of constraints from the formulations (thus avoiding to include them initially), which are later added when necessary in the solution process. This automatic feature of solvers is interesting when the constraint separation must be done by enumeration due to the efficient implementation of the solvers in this case. 

In our particular case,  the use of cut pools in the branch-and-cut method for \eqref{in:3f1} in  formulation $\text{DOMP}_{\text{WOC}}$  seems to require $\mathcal{O}(n^6)$ operations, since there are $\mathcal{O}(n^3)$ \eqref{in:3f1} and $\mathcal{O}(n^3)$ $x$-variables in the formulation. Actually, it is $\mathcal{O}(n^5)$ since each constraint has only $\mathcal{O}(n^2)$ variables to check. Nevertheless, based on the knowledge of the structure of the problem, an efficient separation method of \eqref{in:3f1} constraints can be developed. This \emph{ad hoc}  separation can be included by callbacks allowing a more efficient implementation of the branch-and-cut. 

Focusing in DOMP,  \cite{Labbe2017} propose an algorithm to separate \eqref{in:3f1} with complexity $\mathcal{O}(n^3)$. It is a remarkable quadratic improvement with respect to the pure enumerative approach. Algorithm \ref{algo:separation} shows in detail this separation procedure, based on the calculation of left-hand sides of the possible cuts adding and subtracting two values in each iteration.

In Section \ref{computational}, we will develop a complete analysis of the branch-and-cut method using pools of cuts and the branch-and-cut method using callbacks. We will determine whether or not it is advisable to use a separation algorithm that takes advantage of the knowledge of the problem through callbacks.

\begin{rem}\label{rem1}
	According to previous experiences \citep{Labbe2017, Deleplanque2020}, when valid inequalities \eqref{in:3f1} are embedded in a branch-and-cut procedure over $\text{DOMP}_{\text{WOC}}$ formulation, they should be added at the root node, but not deeper, in order to find a compromise between the integrality gap and the size of the problem. Hence, for our computational study, we will use this  cut-and-branch procedure to check the performance of the solution method proposed in this section.
\end{rem}

\subsection{Row generation procedure for $\text{DOMP}_{\text{SOC}}$}\label{section:rowg}

Since $\text{DOMP}_{\text{WOC}}$ formulation presents coefficients that are not zero-one, in this work we explore the use of  a \eqref{in:3f1} relaxation of $\text{DOMP}_{\text{SOC}}$ adding iteratively these constraints whenever they are necessary. Therefore, the solution method of this section starts with a formulation which has a zero-one coefficient matrix to later add set packing constraints. Hence, we provide a well-behaved (from the solvers point of view) formulation of DOMP without using a huge number of constraints.

The initial formulation which is considered in this row generation procedure is the following.

\begin{align}
\textbf{($\text{DOMP}_{\text{relax}}$)}  \min\quad&\displaystyle\sum_{i\in I}\sum_{j\in I}\sum_{k\in K}\lambda^kc_{ij}x_{ij}^k\nonumber\\
\mbox{s.t.}\quad&\eqref{eq:1fb} -  \eqref{binary:xf1}.\nonumber
\end{align}

Observe that this formulation corresponds to the DOMP model without imposing the order constraints. Therefore, we are dealing with a relaxation of $\text{DOMP}_{\text{SOC}}$. In this case, the proposed relaxation results in the $p$-median problem. 

For each obtained solution in the branch-and-bound of $\text{DOMP}_{\text{relax}}$, \eqref{in:3f1}  are checked by using Algorithm \ref{algo:separation} and added to the model when necessary. Consequently, this row generation method ensures the order by only using a moderate number of \eqref{in:3f1}. This allows to handle bigger instances to be solved in a reasonable computing time as we will see in Section \ref{computational}.

\begin{algorithm}[htbp]
	\textbf{Input:} Let $(\widehat{{\bf x}},\widehat{\bf y})$ be a solution of ($\text{DOMP}_{\text{WOC}}$)/($\text{DOMP}_{\text{relax}}$) with a subset of \eqref{in:3f1}.\\
	\textbf{Output:} Violated cuts / model constraints \eqref{in:3f1}.
	\begin{algorithmic}[1]
		\State Let $(i,j)$ such that $r_{ij}=1$. Then, compute:
		$$lhs:=\sum_{\hat\i\in I}\sum_{\stackrel{\hat\j\in I:}{r_{\hat\i\hat\j}\ge 1}}\widehat{x}_{\hat\i\hat\j}^1+ \widehat{x}_{ij}^{2}.$$
		\If {$lhs>1$}
		\State	Add constraint $\sum_{\hat\i\in I}\sum_{\stackrel{\hat\j\in I:}{r_{\hat\i\hat\j}\ge 1}}x_{\hat\i\hat\j}^1+ x_{ij}^{2}\le 1$.
		\EndIf
		\For{$\ell=2,\ldots,n^2$}
		\State Let $(i,j), (i',j')$ such that $r_{ij}=\ell, r_{i'j'}=\ell-1$. Then, compute:
		$$lhs:=lhs + \widehat{x}_{ij}^{2} - \widehat{x}_{i'j'}^{1}.$$
		\If {$lhs>1$}
		\State	Add constraint $\sum_{\hat\i\in I}\sum_{\stackrel{\hat\j\in I:}{r_{\hat\i\hat\j}\ge\ell}}x_{\hat\i\hat\j}^1+ \sum_{\hat\i\in I}\sum_{\stackrel{\hat\j\in I:}{r_{\hat\i\hat\j}\le\ell}}x_{\hat\i\hat\j}^{2}\le 1$.
		\EndIf
		\EndFor
		\For{$k=2,\ldots,n-1$}
		\State Let $(i,j), (i,'j')$ such that $r_{ij}=1, r_{i'j'}=n^2$. Then, compute:
		$$lhs:=lhs + \widehat{x}_{ij}^{k+1} - \widehat{x}_{i'j'}^{k-1}.$$
		\If {$lhs>1$}
		\State	Add constraint $\sum_{\hat\i\in I}\sum_{\stackrel{\hat\j\in I:}{r_{\hat\i\hat\j}\ge 1}}x_{\hat\i\hat\j}^{k}+ x_{ij}^{k+1}\le 1$.
		\EndIf
		\For{$\ell=2,\ldots,n^2$}
		\State	Let $(i,j), (i',j')$ such that $r_{ij}=\ell, r_{i'j'}=\ell-1$. Then, compute:
		$$lhs:=lhs + \widehat{x}_{ij}^{k+1} - \widehat{x}_{i'j'}^{k}.$$
		\If {$lhs>1$}
		\State	Add constraint $\sum_{\hat\i\in I}\sum_{\stackrel{\hat\j\in I:}{r_{\hat\i\hat\j}\ge\ell}}x_{\hat\i\hat\j}^k+ \sum_{\hat\i\in I}\sum_{\stackrel{\hat\j\in I:}{r_{\hat\i\hat\j}\le\ell}}x_{\hat\i\hat\j}^{k+1}\le 1$.
		\EndIf		
		\EndFor
		\EndFor
	\end{algorithmic}
	{\bf Return:} All violated cuts / model constraints found from \eqref{in:3f1} family.
	
	\caption{(SOC) separation}\label{algo:separation}
\end{algorithm}

\subsubsection{Bounds in constraint relaxations}

In constraint relaxations, any integer solution has to be checked to be valid according to the problem definition. However, there are different alternatives for continuing the branch-and-bound tree exploration when a fractional solution arises. One of them is checking all model constraints to improve the lower bound. Another one is to branch in a particular fractional variable.

One issue to be taken into account is the way in which a subset of \eqref{in:3f1}, that a solution does not verify, is selected to be included in the formulation in order to improve the lower bound without increasing too much the formulation size.
To develop this idea, we introduce the following proposition.

\begin{prop}\label{prop1}
	Given an integer solution $(\widehat{x},\widehat{y})$ for $\text{DOMP}_{\text{relax}}$, if $\widehat{x}$ verifies that
	\begin{equation}\sum_{i'\in I}\sum_{\substack{j'\in I:\\ r_{i'j'}\leq r_{ij}}}x_{i'j'}^k + \sum_{i'\in I}\sum_{\substack{j'\in I:\\ r_{i'j'}\geq r_{ij}}}x_{i'j'}^{k-1}  \le b,\quad i,j \in I,\; k\in K, k \neq 1,\label{socb}
	\end{equation}
	for $b\in [1,2)$, then $\widehat{x}$ satisfies all \eqref{in:3f1}.
\end{prop}
\dem
If $b=1$, then \eqref{socb} are equal to \eqref{in:3f1}. Thus, the result follows trivially. Assume that $b>1$. In this case, since $\widehat{x}$ is integer and $b<2$, then the left hand side of \eqref{socb} must be at most one. Consequently, $\widehat{x}$ satisfies \eqref{in:3f1}.
\fin

As a result of Proposition \ref{prop1}, when an integer solution is obtained in the branch-and-bound tree of $\text{DOMP}_{\text{relax}}$, the $lhs$ described in Algorithm \ref{algo:separation} could be compared with $b\in[1,2)$ instead of comparing it with $1$. Therefore, if $lhs>b$ in Algorithm \ref{algo:separation}, then the corresponding \eqref{in:3f1} are included.

However, varying this $b$ value could affect the number of added cuts when a fractional solution is found and thus, the number of explored nodes in the branch-and bound tree. When $b$ is close to $2$, then the number of identified constraints \eqref{socb} which are not verified by the solution is smaller and they are the most violated cuts. Consequently, for big values of $b$, the number of added cuts will be reduced. Nonetheless, since the number of added cuts is smaller, the number of explored nodes in the branch-and-bound is expected to be bigger.

\begin{rem}
	Following Remark \ref{rem1}, we separate fractional solutions only at root node. Hence, for deeper nodes, Algorithm \ref{algo:separation} is called only when integer solutions are found obtaining upper bounds. Beyond these concerns, lower and upper bounds get closer within the branch-and-bound tree as usually.
\end{rem}

In order to experimentally check how the value of $b$ could impact in times, cuts, and nodes of the row generation proposed in this section, we present Table \ref{tab:accurUB}. We show the results for the instances of sizes $n=20$, $30$, and $40$ that will be detailed in Section \ref{computational}. Particularly, in Table \ref{tab:accurUB}, first column shows the number of clients; the second column reports the number of open facilities; the third set of columns shows the computing time for each $b$ value; the fourth set of columns represents the number of added cuts in the row generation procedure; and finally, the last group of columns reports the number of explored nodes in the branch-and-bound tree. In all cases, each row reports the average value of ten instances. We have tested the results for $b=1$, $b=1.1$, and $b=1.3$.

In Table \ref{tab:accurUB}, we can observe that $b=1$ reports better results than $b=1.1$ and $b=1.3$ when $n=40$. Note that for the largest instances, although the number of added cuts in the branch-and-bound process is bigger, the number of nodes and the computing times are smaller since the gaps at the root node are smaller. Consequently, for the computational results reported in Section \ref{computational}, the choice of $b=1$ is used in the row generation procedure.

\begin{table}
	\caption{Time, cuts, and nodes for different r. h. s. in the row generation procedure using Algorithm \ref{algo:separation} to separate \eqref{socb}\label{tab:accurUB}}}
\centering{
	\begin{adjustbox}{max width=1.0\textwidth}
		\vspace*{-6cm}
		\begin{tabular}{ccrrrrrrrrr}\hline
			& &  \multicolumn{3}{c}{\texttt{Time}} & \multicolumn{3}{c}{\texttt{Cuts}}  & \multicolumn{3}{c}{\texttt{Nodes}} \\
			\cmidrule(lr){3-5}\cmidrule(lr){6-8}\cmidrule(lr){9-11}
			$n$  &$p$ &  \multicolumn{1}{c}{$b=1.0$} &\multicolumn{1}{c}{$b=1.1$} & \multicolumn{1}{c}{$b=1.3$} & \multicolumn{1}{c}{$b=1.0$} &\multicolumn{1}{c}{$b=1.1$} & \multicolumn{1}{c}{$b=1.3$} &  \multicolumn{1}{c}{$b=1.0$} &\multicolumn{1}{c}{$b=1.1$} & \multicolumn{1}{c}{$b=1.3$} \\
			\cmidrule(lr){1-1}
			\cmidrule(lr){2-2}
			\cmidrule(lr){3-3}
			\cmidrule(lr){4-4}
			\cmidrule(lr){5-5}
			\cmidrule(lr){6-6}
			\cmidrule(lr){7-7}
			\cmidrule(lr){8-8}
			\cmidrule(lr){9-9}
			\cmidrule(lr){10-10}
			\cmidrule(lr){11-11}

			\multirow{3}{*}{20}	&	5	&	11.64	&	11.73	&\bf	9.93	&	1440	&	1293	&\bf	1048	&	1	&\bf	1	&	3	\\
			&	6	&\bf	6.47	&	6.83	&	6.82	&	1215	&	1093	&\bf	917	&\bf	1	&	1	&	3	\\
			&	10	&	1.22	&\bf	1.09	&	1.12	&	554	&	512	&\bf	438	&\bf	1	&\bf	1	&\bf	1	\\
			\hline\multirow{3}{*}{30}	&	7	&	158.08	&\bf	137.63	&	506.88	&	4576	&	3930	&\bf	3258	&\bf	11	&	16	&	1133	\\
			&	10	&	75.47	&\bf	69.65	&	131.64	&	2714	&	2412	&\bf	1992	&\bf	77	&	122	&	645	\\
			&	15	&	22.15	&\bf	21.84	&	27.70	&	1429	&	1302	&\bf	1129	&	44	&\bf	43	&	170	\\
			\hline\multirow{3}{*}{40}	&	10	&\bf	576.82	&	859.65	&	4338.31	&	7147	&	6243	&\bf	5308	&\bf	91	&	421	&	5517	\\
			&	13	&\bf	1367.96	&	1564.24	&	4222.45	&	5374	&	4918	&\bf	4064	&\bf	964	&	1268	&	6789	\\
			&	20	&\bf	1060.18	&	1091.21	&	1860.01	&	2691	&	2491	&\bf	2163	&	2408	&\bf	2398	&	5935	\\
			\hline\multicolumn{2}{c}{\bf Total Average:} 			&\bf	328.00	&	376.40	&	1110.52	&	2714	&	2419	&\bf	2032	&\bf	360	&	427	&	2020	\\

			\hline
		\end{tabular}%
	\end{adjustbox}
\end{table}%

\section{Computational experiments} \label{computational}

This section is devoted to the analysis of the solution methods introduced in this paper. The goals of this computational study are the following: 1) checking the differences among the results of formulations  $\text{DOMP}_{\text{SOC}}$ and $\text{DOMP}_{\text{WOC}}$ and determining their limitations; 2) comparing two approaches to implement the branch-and-cut and the row generation algorithms for DOMP.  The difference between these two approaches relies on the fact that the first one uses a fixed constraint/cut pool handled by the solver and the second one applies the separator described in Algorithm \ref{algo:separation} implementing a callback; 3) comparing the different results between the two solution methods defined in Section \ref{section:3}.

The instances used in this computational study were introduced for the first time in \cite{Deleplanque2020}. These instances were created to test different weighted order vectors $\lambda$ beyond the classical ones, namely $p$-median, $p$-center, $k$-centrum problems, etc. The weighted vector $\lambda$ was randomly generated such that $\lambda^k \in \left[\frac{n}{4}, n\right]$ for $ k \in K$. Furthermore, in that data set, there are small- to large-sized instances to perform an exhaustive computational study.  The reader can find the mentioned instances in \citep{DOMPrepository}.

The models were coded in C and solved with SCIP v.6.0.2 \citep{Achterberg,GleixnerEtal2018ZR} using as
optimization solver CPLEX 20.1.0 on a Mac OS Catalina with a Core Intel Xeon W clocked at 3.2 GHz and 96 GB of RAM memory.

In the computational experience, for all the considered formulations and solutions methods, we have included a preprocessing phase. In this stage, we are able to reduce the number of necessary variables to define the problem in terms of optimality. We refer the reader to \cite{Labbe2017} for more details. Besides, we have given an incumbent solution provided by a GRASP heuristic \citep{Deleplanque2020}. This solution let us provide a good upper bound from the beginning of the corresponding solution method. 

Table \ref{tab:formUB} contains the results within two hours of 90 instances up to 40 clients, namely ten instances of each configuration of $n$ and $p$ for two different formulations: $\text{DOMP}_{\text{WOC}}$ and $\text{DOMP}_{\text{SOC}}$. This table and the following ones show the average results for these ten instances: the average CPU time (\texttt{Time}), the number of instances not solved in the time limit (\texttt{\#Unsolved}), the gap at the root node (\texttt{GAProot(\%)}), the gap at termination (\texttt{GAP(\%)}), the number of variables after preprocessing (\texttt{Vars}), the number of constraints (\texttt{OrigCons}), the number of cuts/constraints added in the procedure (\texttt{Cuts}), the number of nodes (\texttt{Nodes}), and the required memory (\texttt{Memory (MB)}). Observe that, for instances with $n=30$ and $n=40$, $\text{DOMP}_{\text{SOC}}$ formulation provides better results than $\text{DOMP}_{\text{WOC}}$. $\text{DOMP}_{\text{SOC}}$ presents a better linear relaxation value (see gap at the root node), and it needs less nodes at the branch-and-bound tree. Consequently,  $\text{DOMP}_{\text{SOC}}$ requires less solution time. In addition, we can conclude that the large number of constraints implies an increment of memory for $\text{DOMP}_{\text{SOC}}$ which makes this formulation too heavy for sizes of $n$ greater than 40. Then, together with the fact that $\text{DOMP}_{\text{WOC}}$ cannot solve any of the instances for $n=40$, we study other alternatives which add (SOC) iteratively until certifying optimality. Thus, in the following tables, we report the results of the methods proposed in Section \ref{section:3}.

\begin{sidewaystable}
	\caption{Results of formulations ($\text{DOMP}_{\text{WOC}}$) and ($\text{DOMP}_{\text{SOC}}$) for $n\in\{ 20,30,40\}$  \label{tab:formUB}}
	\centering{
		\begin{adjustbox}{max width=1.0\textwidth}
			\vspace*{-6cm}
			\begin{tabular}{ccrcrcrrrrrrrrrrrr}\hline
				&& \multicolumn{4}{c}{\texttt{Time (\#Unsolved)}} &  \multicolumn{2}{c}{\texttt{GAProot(\%)}} &  \multicolumn{2}{c}{\texttt{GAP(\%)}}  & \multicolumn{1}{c}{\texttt{Vars}} & \multicolumn{2}{c}{\texttt{OrigCons}}  &  \multicolumn{2}{c}{\texttt{Nodes}}  &  \multicolumn{2}{c}{\texttt{Memory (MB)}}  \\
				\cmidrule(lr){3-6}\cmidrule(lr){7-8}\cmidrule(lr){9-10}\cmidrule(lr){11-11}\cmidrule(lr){12-13}\cmidrule(lr){14-15}\cmidrule(lr){16-17}
				$n$  &$p$ & \multicolumn{2}{c}{$\text{DOMP}_{\text{WOC}}$} & \multicolumn{2}{c}{$\text{DOMP}_{\text{SOC}}$}&  $\text{DOMP}_{\text{WOC}}$ & $\text{DOMP}_{\text{SOC}}$ &  $\text{DOMP}_{\text{WOC}}$ & $\text{DOMP}_{\text{SOC}}$&  &  $\text{DOMP}_{\text{WOC}}$ & $\text{DOMP}_{\text{SOC}}$&  $\text{DOMP}_{\text{WOC}}$ & $\text{DOMP}_{\text{SOC}}$&  $\text{DOMP}_{\text{WOC}}$ & $\text{DOMP}_{\text{SOC}}$\\
				\cmidrule(lr){1-1}
				\cmidrule(lr){2-2}
				\cmidrule(lr){3-4}
				\cmidrule(lr){5-6}
				\cmidrule(lr){7-7}
				\cmidrule(lr){8-8}
				\cmidrule(lr){9-9}
				\cmidrule(lr){10-10}
				\cmidrule(lr){11-11}
				\cmidrule(lr){12-12}
				\cmidrule(lr){13-13}
				\cmidrule(lr){14-14}
				\cmidrule(lr){15-15}
				\cmidrule(lr){16-16}
				\cmidrule(lr){17-17}
				
				\multirow{3}{*}{20}	&	5	&\bf	6.30	&\bf(	0	)&	35.73	&(	0	)&	4.42	&\bf	0.10	&	0.00	&	0.00	&	6054	&\bf	460	&	8041	&	348	&\bf	1	&\bf	39	&	543	\\
				&	6	&\bf	6.10	&\bf(	0	)&	23.19	&(	0	)&	4.37	&\bf	0.01	&	0.00	&	0.00	&	5706	&\bf	460	&	8041	&	297	&\bf	1	&\bf	36	&	518	\\
				&	10	&\bf	1.59	&\bf(	0	)&	9.03	&(	0	)&	2.74	&\bf	0.01	&	0.00	&	0.00	&	4211	&\bf	460	&	8041	&	11	&\bf	1	&\bf	23	&	406	\\
				\hline\multirow{3}{*}{30}	&	7	&	3271.17	&(	2	)&\bf	553.75	&\bf(	0	)&	7.33	&\bf	0.76	&	0.82	&\bf	0.00	&	20643	&\bf	990	&	27061	&	109507	&\bf	6	&\bf	1910	&	3280	\\
				&	10	&	2592.65	&(	1	)&\bf	339.62	&\bf(	0	)&	7.56	&\bf	0.66	&	0.19	&\bf	0.00	&	18245	&\bf	990	&	27061	&	128481	&\bf	16	&\bf	1055	&	2870	\\
				&	15	&	442.39	&(	0	)&\bf	198.13	&\bf(	0	)&	7.43	&\bf	0.30	&	0.00	&	0.00	&	13952	&\bf	990	&	27061	&	22433	&\bf	17	&\bf	261	&	2308	\\
				\hline\multirow{3}{*}{40}	&	10	&	7200.95	&(	10	)&\bf	3437.82	&\bf(	0	)&	9.04	&\bf	1.46	&	6.67	&\bf	0.00	&	48065	&\bf	1720	&	64081	&	56677	&\bf	51	&\bf	5250	&	12724	\\
				&	13	&	7200.87	&(	10	)&\bf	4221.43	&\bf(	3	)&	10.39	&\bf	2.05	&	7.37	&\bf	0.36	&	43664	&\bf	1720	&	64081	&	62211	&\bf	148	&\bf	4335	&	11221	\\
				&	20	&	7200.65	&(	10	)&\bf	3748.44	&\bf(	1	)&	11.62	&\bf	2.46	&	5.45	&\bf	0.42	&	32820	&\bf	1720	&	64081	&	116803	&\bf	331	&\bf	2605	&	8665	\\
				\hline\multicolumn{2}{c}{\bf Total Average:} 			&	3102.52	&(	33	)&\bf	1396.35	&\bf(	4	)&	7.21	&\bf	0.87	&	2.28	&\bf	0.09	&	21484	&\bf	1057	&	33061	&	55196	&\bf	64	&\bf	1724	&	4726	\\

				\hline
			\end{tabular}%
	\end{adjustbox}}
	
\end{sidewaystable}%

Nowadays, commercial solvers include options to code valid inequalities in a branch-and-cut procedure. In this context, these valid inequalites are usually known as \emph{user cuts}, while in constraint relaxations, these model constraints are known as \emph{lazy constraints}. The coding of the cuts/constraints can be done just giving them as input of the linear program. These automatic strategies need to encode all the cuts/constraints in advance within a fixed pool, with an ensuing waste of computer memory. Also, the management of these potential cuts/constraints follows a pre-implemented strategy without taking advantage of the particularities of the formulation beyond the solver pattern recognition based on developers' experience. On the other hand, the use of an oracle (in our case, Algorithm \ref{algo:separation}) to add the cuts/constraints on the fly implementing callbacks could save memory on the whole process \citep[see, e.g.,][and the references therein]{Ackooij2014,Blado2021,Oliveira2014,Mazzi2020,Wolf2014}. Besides, it allows to control when to check and to add those cuts/constraints that is an advantage by itself. For instance, in the row generation solution method, we check model constraints \eqref{in:3f1} at the root node (regardless the solution is fractional or integer) and in any node with integer solution. We refer the reader to \citep{CPLEX,SCIP1,SCIP2} for a detailed discussion.

\begin{sidewaystable}
	\caption{B\&C: $\text{DOMP}_{\text{WOC}}$ with (SOC) as valid inequalities \label{tab:validUB}}
	\centering{
		\begin{adjustbox}{max width=1.0\textwidth}
			\begin{tabular}{ccrcrcrrrrrrrrrrrrr}\hline
				& & \multicolumn{4}{c}{\texttt{Time (\#Unsolved)}} &  \multicolumn{2}{c}{\texttt{GAProot(\%)}} &  \multicolumn{2}{c}{\texttt{GAP(\%)}}  & \multicolumn{2}{c}{\texttt{OrigCons}}& \multicolumn{2}{c}{\texttt{Cuts}}  &  \multicolumn{2}{c}{\texttt{Nodes}}  &  \multicolumn{2}{c}{\texttt{Memory (MB)}}  \\
				\cmidrule(lr){3-6}\cmidrule(lr){7-8}\cmidrule(lr){9-10}\cmidrule(lr){11-12}\cmidrule(lr){13-14}\cmidrule(lr){15-16}\cmidrule(lr){17-18}
				$n$  &$p$ &  \multicolumn{2}{c}{Pool } & \multicolumn{2}{c}{Callback}&  Pool & Callback &  Pool & Callback &  Pool & Callback &  Pool & Callback &  Pool & Callback &  Pool & Callback \\
				\cmidrule(lr){1-1}
				\cmidrule(lr){2-2}
				\cmidrule(lr){3-4}
				\cmidrule(lr){5-6}
				\cmidrule(lr){7-7}
				\cmidrule(lr){8-8}
				\cmidrule(lr){9-9}
				\cmidrule(lr){10-10}
				\cmidrule(lr){11-11}
				\cmidrule(lr){12-12}
				\cmidrule(lr){13-13}
				\cmidrule(lr){14-14}
				\cmidrule(lr){15-15}
				\cmidrule(lr){16-16}
				\cmidrule(lr){17-17}
				\cmidrule(lr){18-18}
				\cmidrule(lr){19-19}
				
				\multirow{3}{*}{20}	&	5	&\bf	15.06	&\bf(	0	)&	15.78	&(	0	)&\bf	0.00	&	0.33	&	0.00	&	0.00	&	8060	&\bf	460	&\bf	544	&	1471	&\bf	1	&	1	&	348	&\bf	106	\\
				&	6	&\bf	9.68	&\bf(	0	)&	9.70	&(	0	)&\bf	0.01	&	0.20	&	0.00	&	0.00	&	8060	&\bf	460	&\bf	485	&	1237	&\bf	1	&	1	&	342	&\bf	101	\\
				&	10	&	3.79	&(	0	)&\bf	1.51	&\bf(	0	)&\bf	0.00	&	0.01	&	0.00	&	0.00	&	8060	&\bf	460	&\bf	297	&	607	&	1	&	1	&	263	&\bf	67	\\
				\hline\multirow{3}{*}{30}	&	7	&	197.98	&(	0	)&\bf	168.44	&\bf(	0	)&\bf	1.02	&	1.08	&	0.00	&	0.00	&	27090	&\bf	990	&\bf	1627	&	4753	&\bf	8	&	20	&	1949	&\bf	430	\\
				&	10	&	103.31	&(	0	)&\bf	66.60	&\bf(	0	)&\bf	1.06	&	1.17	&	0.00	&	0.00	&	27090	&\bf	990	&\bf	1254	&	2953	&\bf	28	&	41	&	1697	&\bf	306	\\
				&	15	&	54.09	&(	0	)&\bf	34.57	&\bf(	0	)&\bf	0.33	&	0.85	&	0.00	&	0.00	&	27090	&\bf	990	&\bf	891	&	2432	&\bf	26	&	67	&	1324	&\bf	266	\\
				\hline\multirow{3}{*}{40}	&	10	&	1279.63	&(	0	)&\bf	626.38	&\bf(	0	)&	4.73	&\bf	1.55	&	0.00	&	0.00	&	64120	&\bf	1720	&\bf	4300	&	7899	&	153	&\bf	95	&	7142	&\bf	1041	\\
				&	13	&	2316.83	&(	1	)&\bf	1343.62	&\bf(	0	)&	5.43	&\bf	2.22	&	0.21	&\bf	0.00	&	64120	&\bf	1720	&\bf	5140	&	6445	&\bf	688	&	890	&	6410	&\bf	861	\\
				&	20	&	1418.48	&(	1	)&\bf	717.14	&\bf(	0	)&	3.82	&\bf	2.78	&	0.24	&\bf	0.00	&	64120	&\bf	1720	&\bf	2839	&	4636	&\bf	813	&	1126	&	4897	&\bf	641	\\
				\hline\multirow{3}{*}{50}	&	12	&	3223.91	&(	0	)&\bf	1858.95	&\bf(	0	)&	2.12	&\bf	0.89	&	0.00	&	0.00	&	125150	&\bf	2650	&\bf	4393	&	12466	&\bf	90	&	150	&	21412	&\bf	2228	\\
				&	16	&	4057.82	&(	2	)&\bf	2364.00	&\bf(	1	)&	3.41	&\bf	1.11	&	0.17	&\bf	0.05	&	125150	&\bf	2650	&\bf	4756	&	10172	&\bf	341	&	658	&	19636	&\bf	1922	\\
				&	25	&\bf	3166.70	&\bf(	1	)&	2311.38	&(	2	)&	4.01	&\bf	1.27	&	0.40	&\bf	0.35	&	125150	&\bf	2650	&\bf	3273	&	9787	&\bf	611	&	912	&	14404	&\bf	1512	\\
				\hline\multirow{3}{*}{60}	&	15	&	6953.04	&(	8	)&\bf	5884.18	&\bf(	5	)&	5.64	&\bf	1.14	&	3.43	&\bf	0.71	&	216180	&\bf	3780	&\bf	4239	&	23249	&	62	&\bf	58	&	51404	&\bf	5013	\\
				&	20	&	6762.80	&(	7	)&\bf	5267.87	&\bf(	6	)&	5.89	&\bf	1.06	&	2.77	&\bf	0.53	&	216180	&\bf	3780	&\bf	4989	&	22206	&\bf	111	&	144	&	46079	&\bf	4528	\\
				&	30	&	6950.37	&(	9	)&\bf	5749.35	&\bf(	6	)&	7.61	&\bf	1.74	&	1.88	&\bf	0.72	&	216180	&\bf	3780	&\bf	5655	&	15535	&\bf	261	&	1169	&	35170	&\bf	3011	\\
				\hline\multicolumn{2}{c}{\bf Total Average:} 			&	2434.23	&(	29	)&\bf	1761.30	&\bf(	20	)&	3.01	&\bf	1.16	&	0.61	&\bf	0.16	&	88120	&\bf	1920	&\bf	2979	&	8390	&\bf	213	&	356	&	14165	&\bf	1469	\\
				\hline
			\end{tabular}%
		\end{adjustbox}
	}
\end{sidewaystable}%

\begin{sidewaystable}
	\caption{Row generation, i.e, SOC (model constraints) added iteratively\label{tab:feasUB}}
	\centering{
		\begin{adjustbox}{max width=1.0\textwidth}
			\vspace*{-6cm}
			\begin{tabular}{ccrcrcrrrrrrrrrrrrr}\hline
				& & \multicolumn{4}{c}{\texttt{Time (\#Unsolved)}} &  \multicolumn{2}{c}{\texttt{GAProot(\%)}} &  \multicolumn{2}{c}{\texttt{GAP(\%)}}  & \multicolumn{2}{c}{\texttt{OrigCons}}& \multicolumn{2}{c}{\texttt{Cuts}}  &  \multicolumn{2}{c}{\texttt{Nodes}}  &  \multicolumn{2}{c}{\texttt{Memory (MB)}}  \\
				\cmidrule(lr){3-6}\cmidrule(lr){7-8}\cmidrule(lr){9-10}\cmidrule(lr){11-12}\cmidrule(lr){13-14}\cmidrule(lr){15-16}\cmidrule(lr){17-18}
				$n$  &$p$ &  \multicolumn{2}{c}{Pool } & \multicolumn{2}{c}{Callback}&  Pool & Callback &  Pool & Callback &  Pool & Callback &  Pool & Callback &  Pool & Callback &  Pool & Callback \\
				\cmidrule(lr){1-1}
				\cmidrule(lr){2-2}
				\cmidrule(lr){3-4}
				\cmidrule(lr){5-6}
				\cmidrule(lr){7-7}
				\cmidrule(lr){8-8}
				\cmidrule(lr){9-9}
				\cmidrule(lr){10-10}
				\cmidrule(lr){11-11}
				\cmidrule(lr){12-12}
				\cmidrule(lr){13-13}
				\cmidrule(lr){14-14}
				\cmidrule(lr){15-15}
				\cmidrule(lr){16-16}
				\cmidrule(lr){17-17}
				\cmidrule(lr){18-18}
				\cmidrule(lr){19-19}
				
				\multirow{3}{*}{20}	&	5	&	15.63	&(	0	)&\bf	11.64	&\bf(	0	)&\bf	0.19	&	0.29	&	0.00	&	0.00	&	8041	&\bf	441	&\bf	514	&	1440	&\bf	1	&	1	&	346	&\bf	129	\\
				&	6	&	10.28	&(	0	)&\bf	6.47	&\bf(	0	)&	0.03	&\bf	0.00	&	0.00	&	0.00	&	8041	&\bf	441	&\bf	483	&	1215	&\bf	1	&	1	&	340	&\bf	124	\\
				&	10	&	3.80	&(	0	)&\bf	1.22	&\bf(	0	)&\bf	0.00	&	0.01	&	0.00	&	0.00	&	8041	&\bf	441	&\bf	309	&	554	&	1	&	1	&	260	&\bf	91	\\
				\hline\multirow{3}{*}{30}	&	7	&\bf	150.01	&\bf(	0	)&	158.08	&(	0	)&\bf	0.69	&	1.11	&	0.00	&	0.00	&	27061	&\bf	961	&\bf	1422	&	4576	&\bf	4	&	11	&	1928	&\bf	497	\\
				&	10	&	88.03	&(	0	)&\bf	75.47	&\bf(	0	)&\bf	0.96	&	1.15	&	0.00	&	0.00	&	27061	&\bf	961	&\bf	1082	&	2714	&\bf	15	&	77	&	1668	&\bf	395	\\
				&	15	&	49.04	&(	0	)&\bf	22.15	&\bf(	0	)&\bf	0.42	&	0.85	&	0.00	&	0.00	&	27061	&\bf	961	&\bf	783	&	1429	&\bf	12	&	44	&	1321	&\bf	367	\\
				\hline\multirow{3}{*}{40}	&	10	&	837.75	&(	0	)&\bf	576.82	&\bf(	0	)&\bf	1.51	&	1.57	&	0.00	&	0.00	&	64081	&\bf	1681	&\bf	2907	&	7147	&\bf	44	&	91	&	7059	&\bf	1256	\\
				&	13	&\bf	1061.32	&\bf(	0	)&	1367.96	&(	1	)&\bf	2.04	&	2.23	&\bf	0.00	&	0.07	&	64081	&\bf	1681	&\bf	3212	&	5374	&\bf	207	&	964	&	6331	&\bf	1099	\\
				&	20	&\bf	1093.46	&\bf(	0	)&	1060.18	&(	1	)&\bf	2.55	&	2.81	&\bf	0.00	&	0.13	&	64081	&\bf	1681	&\bf	2404	&	2691	&\bf	680	&	2408	&	4839	&\bf	877	\\
				\hline\multirow{3}{*}{50}	&	12	&	3061.61	&(	1	)&\bf	1648.54	&\bf(	0	)&	2.08	&\bf	0.88	&	0.12	&\bf	0.00	&	125101	&\bf	2601	&\bf	3858	&	10720.2	&\bf	91.7	&	116	&	21271	&\bf	2609	\\
				&	16	&	3104.23	&(	0	)&\bf	1691.65	&\bf(	0	)&\bf	1.09	&	1.11	&	0.00	&	0.00	&	125101	&\bf	2601	&\bf	3429	&	7236.6	&\bf	219.9	&	662	&	19496	&\bf	2248	\\
				&	25	&	2958.29	&(	2	)&\bf	1838.77	&\bf(	1	)&\bf	1.23	&	1.26	&	0.34	&\bf	0.28	&	125101	&\bf	2601	&\bf	2993	&	3919.9	&\bf	526	&	2217	&	14316	&\bf	1752	\\
				\hline\multirow{3}{*}{60} 	&	15	&	7313.70	&(	9	)&\bf	5357.93	&\bf(	5	)&	4.45	&\bf	1.15	&	3.19	&\bf	0.67	&	216121	&\bf	3721	&\bf	4515	&	17078	&\bf	9.2	&	63	&	51469	&\bf	5346	\\
				&	20	&	6439.44	&(	6	)&\bf	4529.22	&\bf(	4	)&	1.45	&\bf	1.05	&	0.94	&\bf	0.36	&	216121	&\bf	3721	&\bf	4605	&	11582	&\bf	33.6	&	444	&	46014	&\bf	4530	\\
				&	30	&	6733.81	&(	8	)&\bf	5093.43	&\bf(	4	)&	2.42	&\bf	1.74	&	1.02	&\bf	0.46	&	216121	&\bf	3721	&\bf	4402	&	6169.1	&\bf	256.6	&	2127	&	34981	&\bf	3463	\\
				\hline\multicolumn{2}{c}{\bf Total Average:} 			&	2194.69	&(	26	)&\bf	1562.64	&\bf(	16	)&	1.41	&\bf	1.15	&	0.37	&\bf	0.13	&	88081	&\bf	1881	&\bf	2461	&	5590	&\bf	140	&	615	&	14109	&\bf	1652	\\
				\hline
			\end{tabular}%
		\end{adjustbox}
	}
\end{sidewaystable}%

Table \ref{tab:validUB} presents the results, in two hours of time limit, of the branch-and-cut introduced in Section \ref{section:B&C}, i.e., ($\text{DOMP}_{\text{WOC}}$) with (SOC) as valid inequalities. Here, we follow two different strategies: we code (SOC) defining a fixed user cut pool in SCIP and the solver decides when to check and to add them (Pool); or we check the constraints using Algorithm \ref{algo:separation} which adds them by an user callback when needed at the root node (Callback). The results exhibit that the method using Algorithm \ref{algo:separation} shows better solving times than the automatic approach. Besides, the method using the callback can solve nine more instances than the automatic one. Note that these better results are explained by the efficiency of  our separation algorithm. Regarding the required memory, observe that memory used by the automatic method increases quickly with $n$.  Therefore, the application of this method does not seem to be useful for bigger instances since they could not be loaded: it requires around 50 GB of RAM memory already for $n=60$.

Table \ref{tab:feasUB} reports the results, in two hours of computing time, of the row generation method, i.e., $\text{DOMP}_{\text{relax}}$ with \eqref{in:3f1}   as model constraints not included from the beginning. Two approaches to carry out the row generation are considered: the automatic use of \eqref{in:3f1} defining a fixed lazy constraint pool (Pool) and the application of Algorithm \ref{algo:separation} to add \eqref{in:3f1} when necessary (Callback). In this table, we have included the same columns as in Table \ref{tab:validUB}. Note that the callback approach provides the best computing time results and only 16 instances remain unsolved after the time limit. The automatic method requires less cuts and nodes than Callback. However, the  required memory increases faster when using the automatic approach because it needs to encode all the original \eqref{in:3f1} constraints which are $\mathcal{O}(n^3)$. Consequently, the performance of the row generation following the automatic separation shows, in general, worse results than using Algorithm \ref{algo:separation}.

From tables \ref{tab:validUB} and \ref{tab:feasUB}, we can conclude that the performance of the automatic branch-and-cut  and the automatic row generation are quite limited in comparison with the use of the separation presented in Algorithm \ref{algo:separation}. The reason of the better performance of the callback approach is that Algorithm \ref{algo:separation} exploits the knowledge of the problem. Therefore, for the study of larger instances, we focus on the branch-and-cut  and row generation approaches using Algorithm \ref{algo:separation}.

Table \ref{tab:5hours} shows a comparison between the results provided by the callback versions of the branch-and-cut method (B\&C) and the row generation technique (RowGen). For these experiments, we establish a time limit of five hours. Overall, the row generation technique outperforms the branch-and-cut in terms of computational times. Moreover,  the row generation is able to solve more instances than the  branch-and-cut  in five hours. 

\begin{sidewaystable}
	\caption{B\&C vs row generation in 5 hours of time limit\label{tab:5hours}}
	\centering{
		\begin{adjustbox}{max width=1.0\textwidth}
			\vspace*{-6cm}
			\begin{tabular}{ccrcrcrrrrrrrrrrrrrr}\hline
				& & \multicolumn{4}{c}{\texttt{Time (\#Unsolved)}} &  \multicolumn{2}{c}{\texttt{GAProot(\%)}} &  \multicolumn{2}{c}{\texttt{GAP(\%)}}  & \multicolumn{1}{c}{\texttt{Vars}} & \multicolumn{2}{c}{\texttt{OrigCons}}& \multicolumn{2}{c}{\texttt{Cuts}}  &  \multicolumn{2}{c}{\texttt{Nodes}}  &  \multicolumn{2}{c}{\texttt{Memory (MB)}}  \\
				\cmidrule(lr){3-6}\cmidrule(lr){7-8}\cmidrule(lr){9-10}\cmidrule(lr){11-11}\cmidrule(lr){12-13}\cmidrule(lr){14-15}\cmidrule(lr){16-17}\cmidrule(lr){18-19}
				$n$  &$p$&  \multicolumn{2}{c}{B\&C} & \multicolumn{2}{c}{RowGen}&  B\&C & RowGen &  B\&C & RowGen&&  B\&C & RowGen&  B\&C & RowGen&  B\&C & RowGen&  B\&C & RowGen\\
				\cmidrule(lr){1-1}
				\cmidrule(lr){2-2}
				\cmidrule(lr){3-4}
				\cmidrule(lr){5-6}
				\cmidrule(lr){7-7}
				\cmidrule(lr){8-8}
				\cmidrule(lr){9-9}
				\cmidrule(lr){10-10}
				\cmidrule(lr){11-11}
				\cmidrule(lr){12-12}
				\cmidrule(lr){13-13}
				\cmidrule(lr){14-14}
				\cmidrule(lr){15-15}
				\cmidrule(lr){16-16}
				\cmidrule(lr){17-17}
				\cmidrule(lr){18-18}
				\cmidrule(lr){19-19}
				
				\multirow{3}{*}{20}	&	5	&	15.78	&(	0	)&\bf	11.64	&\bf(	0	)&	0.33	&\bf	0.29	&	0.00	&	0.00	&	6054	&	460	&\bf	441	&	1471	&\bf	1440	&\bf	1	&	1	&\bf	106	&	129	\\
				&	6	&	9.70	&(	0	)&\bf	6.47	&\bf(	0	)&	0.20	&\bf	0.00	&	0.00	&	0.00	&	5706	&	460	&\bf	441	&	1237	&\bf	1215	&	1	&	1	&\bf	101	&	124	\\
				&	10	&	1.51	&(	0	)&\bf	1.22	&\bf(	0	)&	0.01	&\bf	0.01	&	0.00	&	0.00	&	4211	&	460	&\bf	441	&	607	&\bf	554	&	1	&	1	&\bf	67	&	91	\\
				\hline\multirow{3}{*}{30}	&	7	&	168.44	&(	0	)&\bf	158.08	&\bf(	0	)&\bf	1.08	&	1.11	&	0.00	&	0.00	&	20643	&	990	&\bf	961	&	4753	&\bf	4576	&	20	&\bf	11	&\bf	430	&	497	\\
				&	10	&\bf	66.60	&\bf(	0	)&	75.47	&(	0	)&	1.17	&\bf	1.15	&	0.00	&	0.00	&	18245	&	990	&\bf	961	&	2953	&\bf	2714	&\bf	41	&	77	&\bf	306	&	395	\\
				&	15	&	34.57	&(	0	)&\bf	22.15	&\bf(	0	)&\bf	0.85	&	0.85	&	0.00	&	0.00	&	13952	&	990	&\bf	961	&	2432	&\bf	1429	&	67	&\bf	44	&\bf	266	&	367	\\
				\hline\multirow{3}{*}{40}	&	10	&	626.38	&(	0	)&\bf	576.82	&\bf(	0	)&\bf	1.55	&	1.57	&	0.00	&	0.00	&	48065	&	1720	&\bf	1681	&	7899	&\bf	7147	&	95	&\bf	91	&\bf	1041	&	1256	\\
				&	13	&\bf	1343.62	&\bf(	0	)&	1390.00	&(	0	)&\bf	2.22	&	2.23	&	0.00	&	0.00	&	43664	&	1720	&\bf	1681	&	6445	&\bf	5374	&\bf	890	&	981	&\bf	861	&	1099	\\
				&	20	&\bf	717.14	&\bf(	0	)&	1082.85	&(	0	)&\bf	2.78	&	2.81	&	0.00	&	0.00	&	32820	&	1720	&\bf	1681	&	4636	&\bf	2691	&\bf	1126	&	2546	&\bf	641	&	877	\\
				\hline\multirow{3}{*}{50}	&	12	&	1858.95	&(	0	)&\bf	1648.54	&\bf(	0	)&	0.89	&\bf	0.88	&	0.00	&	0.00	&	94784	&	2650	&\bf	2601	&	12466	&\bf	10720	&	150	&\bf	116	&\bf	2228	&	2609	\\
				&	16	&	2497.71	&(	0	)&\bf	1691.65	&\bf(	0	)&	1.11	&\bf	1.11	&	0.00	&	0.00	&	85630	&	2650	&\bf	2601	&	10172	&\bf	7237	&	687	&\bf	662	&\bf	1923	&	2248	\\
				&	25	&	3425.89	&(	1	)&\bf	2918.77	&\bf(	1	)&	1.27	&\bf	1.26	&	0.29	&\bf	0.23	&	63776	&	2650	&\bf	2601	&	9787	&\bf	3920	&\bf	1103	&	3277	&\bf	1516	&	1777	\\
				\hline\multirow{3}{*}{60} 	&	15	&	11007.79	&(	4	)&\bf	9004.34	&\bf(	2	)&\bf	1.13	&	1.14	&	0.33	&\bf	0.30	&	161807	&	3780	&\bf	3721	&	23249	&\bf	17095	&	291	&\bf	268	&\bf	5046	&	5438	\\
				&	20	&	9238.02	&(	3	)&\bf	6695.45	&\bf(	1	)&	1.06	&\bf	1.05	&	0.28	&\bf	0.17	&	144983	&	3780	&\bf	3721	&	22206	&\bf	11582	&\bf	353	&	599	&	4542	&\bf	4540	\\
				&	30	&	10103.21	&(	3	)&\bf	8524.74	&\bf(	3	)&	1.74	&\bf	1.74	&	0.39	&\bf	0.24	&	109804	&	3780	&\bf	3721	&	15535	&\bf	6169	&\bf	2071	&	3381	&\bf	3030	&	3493	\\
				\hline\multirow{3}{*}{70} 	&	17	&	16541.03	&(	8	)&\bf	15920.52	&\bf(	8	)&	1.01	&\bf	0.95	&	0.67	&\bf	0.56	&	259406	&	5110	&\bf	5041	&	31406	&\bf	23228	&\bf	95	&	161	&\bf	9281	&	9666	\\
				&	23	&	16506.34	&(	8	)&\bf	15748.13	&\bf(	8	)&\bf	1.10	&	1.10	&	0.68	&\bf	0.54	&	231680	&	5110	&\bf	5041	&	33833	&\bf	16004	&\bf	215	&	654	&	8792	&\bf	7851	\\
				&	35	&	18011.10	&(	10	)&\bf	17281.08	&\bf(	8	)&	2.05	&\bf	2.04	&	1.34	&\bf	1.13	&	173955	&	5110	&\bf	5041	&	22814	&\bf	8648	&\bf	721	&	1924	&\bf	5686	&	6259	\\
				\hline\multirow{3}{*}{80} 	&	20	&	18041.98	&(	10	)&\bf	18001.87	&\bf(	10	)&	4.22	&\bf	1.47	&	4.22	&\bf	1.37	&	383199	&	6640	&\bf	6561	&	42109	&\bf	30421	&\bf	1	&	11	&	14227	&\bf	13641	\\
				&	26	&	18030.43	&(	10	)&\bf	17433.69	&\bf(	9	)&	1.39	&\bf	0.78	&	1.32	&\bf	0.57	&	346926	&	6640	&\bf	6561	&	38792	&\bf	20784	&\bf	100	&	325	&	13194	&\bf	11930	\\
				&	40	&	16047.37	&(	7	)&\bf	14613.12	&\bf(	6	)&	0.57	&\bf	0.57	&	0.32	&\bf	0.29	&	259186	&	6640	&\bf	6561	&	17419	&\bf	10420	&\bf	794	&	1273	&\bf	6530	&	8940	\\
				\hline\multirow{3}{*}{90} 	&	22	&	18103.76	&(	10	)&\bf	18002.56	&\bf(	10	)&	8.21	&\bf	7.34	&	8.21	&\bf	7.34	&	549561	&	8370	&\bf	8281	&	28820	&\bf	28267	&	1	&	1	&	14849	&\bf	13946	\\
				&	30	&	18078.18	&(	10	)&\bf	17886.15	&\bf(	9	)&	4.35	&\bf	0.56	&	4.34	&\bf	0.46	&	488316	&	8370	&\bf	8281	&	35203	&\bf	23544	&\bf	3	&	37	&\bf	15053	&	16879	\\
				&	45	&	14629.97	&(	7	)&\bf	14340.26	&\bf(	6	)&\bf	0.56	&	0.56	&	0.36	&\bf	0.31	&	368560	&	8370	&\bf	8281	&	18634	&\bf	13363	&\bf	410	&	700	&\bf	8715	&	13660	\\
				\hline\multirow{3}{*}{100} 	&	25	&	18150.10	&(	10	)&\bf	18006.01	&\bf(	10	)&\bf	8.01	&	11.71	&\bf	8.01	&	11.71	&	749074	&	10300	&\bf	10201	&\bf	19420	&	21398	&	1	&	1	&	14407	&\bf	13231	\\
				&	33	&	18138.18	&(	10	)&\bf	18004.78	&\bf(	10	)&	7.05	&\bf	3.88	&	7.05	&\bf	3.88	&	672511	&	10300	&\bf	10201	&	24045	&\bf	23931	&	1	&	1	&\bf	14681	&	16009	\\
				&	50	&	18065.05	&(	10	)&\bf	17980.43	&\bf(	9	)&	0.58	&\bf	0.58	&	0.53	&\bf	0.44	&	504981	&	10300	&\bf	10201	&	25922	&\bf	16253	&\bf	66	&	351	&\bf	13819	&	18941	\\
				\hline\multicolumn{2}{c}{\bf Total Average:} 			&	9239.22	&(	121	)&\bf	8778.77	&\bf(	110	)&	2.09	&\bf	1.80	&	1.42	&\bf	1.09	&	216352	&	4447	&\bf	4388	&	17195	&\bf	11856	&\bf	345	&	648	&\bf	5975	&	6515	\\
				\hline
			\end{tabular}%
	\end{adjustbox}}
	
\end{sidewaystable}%

\begin{sidewaystable}
	\caption{B\&P\&C, B\&C, and row generation in 24 hours of time limit\label{tab:24hours}}
	\centering{
		\begin{adjustbox}{max width=1.0\textwidth}
			\vspace*{-6cm}
			\begin{tabular}{clrrrrrrrrrrrrrr}
				\toprule
				\multicolumn{2}{c}{$n$} & \multicolumn{6}{c}{100}                     & \multicolumn{6}{c}{120}                     &  \\
				\cmidrule(lr){1-2}
				\cmidrule(lr){3-8}
				\cmidrule(lr){9-14}
				\multicolumn{2}{c}{$p$} & \multicolumn{2}{c}{25} & \multicolumn{2}{c}{33} & \multicolumn{2}{c}{50} & \multicolumn{2}{c}{30} & \multicolumn{2}{c}{40} & \multicolumn{2}{c}{60} & \multicolumn{2}{c}{\bf Total Average}\\
				\cmidrule(lr){1-2}
				\cmidrule(lr){3-4}
				\cmidrule(lr){5-6}
				\cmidrule(lr){7-8}
				\cmidrule(lr){9-10}
				\cmidrule(lr){11-12}
				\cmidrule(lr){13-14}
				\cmidrule(lr){15-16}
				\multicolumn{1}{c}{\multirow{3}[2]{*}{\texttt{Time (\#Unsolved)}}} & B\&P\&C &    	86400.19& (10) &	86400.29& (10) &	86400.08& (10) &86400.18& (10) &	86400.27& (10) &	86400.12& (10) &	86400.19& (60)
				\\
				& B\&C  & 86557.06 & (10)  & 85594.56 & (9)   & \textbf{78019.88} & \textbf{(8)} & 86841.70 & (10)  & 86833.42 & (10)  & 81439.12 & (9)   & 84214.29 & (56) \\
				& RowGen & \textbf{74555.59} & \textbf{(7)} & \textbf{76691.23} & \textbf{(8)} & 79403.08 & (8)   & \textbf{86403.24} & \textbf{(10)} & \textbf{86402.13} & \textbf{(10)} & \textbf{79855.20} & \textbf{(9)} & \textbf{80551.74} & \textbf{(52)} \\
				\midrule
				\multicolumn{1}{c}{\multirow{3}[2]{*}{\texttt{GAProot (\%)}}} & B\&P\&C&  4.50& &	3.65& &	2.32& &	\textbf{5.86}& &4.59& &2.50& &3.90& 
				\\
				& B\&C  & 4.84 &         & 1.08  &       & 0.58  &       & 8.74  &       & 7.51  &       & \textbf{0.49} &       & 3.87  &  \\
				& RowGen & \textbf{0.57} &       & \textbf{0.66} &       & \textbf{0.50} &       & 7.60 &       & \textbf{2.06} &       & 0.49  &       & \textbf{1.98} &  \\
				\midrule
				\multicolumn{1}{c}{\multirow{3}[2]{*}{\texttt{GAP (\%)}}} & B\&P\&C & 4.50& &3.65& &	2.32& &	\textbf{5.86}& &4.59& &	2.50& &	3.90& 
				\\
				& B\&C  & 4.83  &       & 1.01  &       & 0.40  &       & 8.74  &       & 7.51  &       & 0.40  &       & 3.81  &  \\
				& RowGen & \textbf{0.36} &       & \textbf{0.52} &       & \textbf{0.28} &       & 7.60 &       & \textbf{2.04} &       & \textbf{0.37} &       & \textbf{1.86} &  \\
				\midrule
				\multicolumn{1}{c}{\multirow{3}[2]{*}{\texttt{Vars} }} & B\&P\&C & \textbf{	47921}& &\textbf{	45038}& &\textbf{	40436}& &\textbf{	52678}& &\textbf{	49576}& &\textbf{	41743}& &\textbf{	46232}& 
				
				\\
				& B\&C  & 595107 &       & 725934 &       & 605525 &       & 1293335 &       & 1155202 &       & 871156 &       & 874377 &  \\
				& RowGen & 595107 &       & 725934 &       & 605525 &       & 1293335 &       & 1155202 &       & 871156 &       & 874377 &  \\
				\midrule
				\multicolumn{1}{c}{\multirow{3}[2]{*}{\texttt{OrigCons }}} &B\&P\&C & \textbf{	301}& &\textbf{	301}& &\textbf{	301}& &\textbf{	361}& &\textbf{	361}& &\textbf{	361}& &\textbf{	331}& 
				
				\\
				& B\&C  & 10300 &       & 10300 &       & 10300 &       & 14760 &       & 14760 &       & 14760 &       & 12530 &  \\
				& RowGen & 10201 &       & 10201 &       & 10201 &       & 14641 &       & 14641 &       & 14641 &       & 12421 &  \\
				\midrule
				\multicolumn{1}{c}{\multirow{3}[2]{*}{\texttt{Cuts} }} & B\&P\&C &     \textbf{	10573}& &\textbf{	9214}& &\textbf{	5760}& &\textbf{	13717}& &\textbf{	12259}& &\textbf{	7611}& &\textbf{	9856}& 
				\\
				& B\&C  & 89913 &       & 77195 &       & 25922 &       & 44518 &       & 55814 &       & 34981 &       & 54724 &  \\
				& RowGen & 25990 &       & 49665 &       & 26137 &       & 44244 &       & 42998 &       & 24038 &       & 35512 &  \\
				\midrule
				\multicolumn{1}{c}{\multirow{3}[2]{*}{\texttt{Nodes} }} & B\&P\&C &      \textbf{	1}& &\textbf{	1}& &\textbf{	1}& &\textbf{	1}& &\textbf{	1}& &\textbf{	1}& &\textbf{	1}& 
				\\
				& B\&C  & 2 &       &58 &       &1205 &       & \textbf{	1}     &       &\textbf{	1} &       & 140 &       & 235 &  \\
				& RowGen & 1198  &       & 88    &       & 1245  &       & \textbf{	1}     &       & 26    &       & 324   &       & 480   &  \\
				\midrule
				\multicolumn{1}{c}{\multirow{3}[2]{*}{\texttt{Memory (MB)} }} & B\&P\&C &      \textbf{	1739}& &\textbf{	1215}& &\textbf{	663}& &\textbf{	2105}& &\textbf{	1413}& &\textbf{	697}& &\textbf{	1305}& 
				\\
				& B\&C  & 46512 &       & 37105 &       & 13992 &       & 39891 &       & 40861 &       & 25062 &       & 33904 &  \\
				& RowGen & 23225 &       & 32057 &       & 23869 &       & 35066 &       & 40005 &       & 35755 &       & 31663 &  \\
				\bottomrule
			\end{tabular}%
	\end{adjustbox}}
\end{sidewaystable}%

Note that, for some huge instances, the problem does not get even the root node bounds. For those instances, $\texttt{GAProot(\%)}$ and $\texttt{GAP(\%)}$ report the same value and thus, $\texttt{GAProot(\%)}$ cannot be analyzed as the gap at the root node, but as the gap at the root node at the time limit.

Regarding $n=60$ instances, the branch-and-cut is able to solve 20 out of 30 instances in five hours, whereas this solution method certifies optimality for only 13 instances  in two hours (see Table \ref{tab:validUB}). These three extra hours also let the row generation algorithm to solve 24 instances, seven more than the same algorithm in two hours (see Table \ref{tab:feasUB}).

For most of the instances, the integrality gap at termination provided by the row generation procedure is smaller than the one obtained by the branch-and-cut, even with less cuts added. This gives us the idea that the added cuts are more accurate when \eqref{in:3f2} family is not included in the formulation. However, for $n=100$, the gap at termination is smaller for the branch-and-cut  procedure since the added cuts cannot improve the lower bound given by the linear relaxation of the program in both solution methods within the time limit.

To analyze the differences between the branch-and-cut and the row generation, starting from $\text{DOMP}_{\text{WOC}}$ and $\text{DOMP}_{\text{relax}}$, respectively, we give the solver up to 24 hours of time limit. Thereby, the influence of the linear relaxation bound is not so decisive. Furthermore, the branch-and-price-and-cut (B\&P\&C) described in  \cite{Deleplanque2020} has been tested for those instances in the same computer and with the same time limit. The reader could see, in Table \ref{tab:24hours}, how the gaps at termination are smaller on average for the row generation procedure. In fact, they are reduced by half and for the  row generation are less than 2\%. This approach needs less cuts to have a reasonable bound at the root node what let it branch faster to generate more nodes and improve the bounds. Thus, whereas the solution method detailed in Section \ref{section:B&C} solves four instances out of 60 for these large-sized instances, the solution method proposed in Section \ref{section:rowg} is able to solve eight instances to optimality and on top of that, the row generation procedure also reduces the gap at termination for those instances which are not solved. The B\&P\&C procedure uses less variables but the gap at termination is worse because it is not able to solve even the root node. However, this approach would still be useful for bigger instances since it requires much less memory.

\section{Conclusions}\label{section:5}
In this work, we have introduced a row generation solution method which has improved the best known performances for DOMP regarding the medium-sized instances of the data set described in \cite{Deleplanque2020}. In adittion, we have improved the best known solution for an instance with $n=90$ and $p=45$ ({\tt domp90p45v5.domp}) that can be found in the mentioned dataset.

For large-sized instances, the lower bound provided by formulations which include \eqref{in:3f2} makes them  also a good alternative. For these instances, the lower bound given by the linear relaxation can be barely improved  within the time limit. Moreover, comparing with the integrality gap given by the branch-and-price algorithm \citep{Deleplanque2020}, one should note that the column generation of its master problem gives theoretically better lower bounds.

Taking into account that the limits of our row generation algorithm come from the huge number of variables, a combination of row and column generation seems to be a promising approach to be considered as future research line.

\section*{Acknowledgments}

The authors of this research acknowledge financial support by the Spanish Ministerio de Ciencia y Tecnolog\'ia, Agencia Estatal de Investigaci\'on and Fondos Europeos de Desarrollo Regional (FEDER) via projects PID2020-114594GB-C21/C22. The authors also acknowledge partial support from NetmeetData: Ayudas Fundación BBVA a equipos de investigación cient\'ifica
2019, Junta de Andalucía P18-FR-1422, and project CIGE/2021/161 from Generalitat Valenciana. The first author thanks partial support from project FEDER-UCA18-106895.  The second and third authors also acknowledge partial support from projects PID2020-114594GB-C21, B-FQM-322-UGR20, FEDER-US-1256951, and FQM-331. The first author has been supported
by Subprograma Juan de la Cierva Formación 2019 (FJC2019-
039023-I) and PAIDI 2020, postdoctoral fellowship financed by
European Social Fund and Junta de Andalucía. The second author also acknowledges PAIDI 2020, postdoctoral fellowship financed by
European Social Fund and Junta de Andalucía. 

We would like to thank Reviewers for taking the necessary time and effort to make constructive suggestions that help to improve the quality of the manuscript.


	

\begin{thebibliography}{33}
	\providecommand{\natexlab}[1]{#1}
	\providecommand{\url}[1]{\texttt{#1}}
	\expandafter\ifx\csname urlstyle\endcsname\relax
	\providecommand{\doi}[1]{doi: #1}\else
	\providecommand{\doi}{doi: \begingroup \urlstyle{rm}\Url}\fi
	
	\bibitem[Achterberg(2009)]{Achterberg}
	Achterberg~T (2009)
	\newblock SCIP: Solving constraint integer programs.
	\newblock Mathematical Programming Computation 1(1):1--41. \url{https://doi.org/10.1007/s12532-008-0001-1}
	
	\bibitem[Ackooij and de~Oliveira(2014)]{Ackooij2014}
	Ackooij~W, de~Oliveira~W (2014)
	\newblock Level bundle methods for constrained convex optimization with various
	oracles.
	\newblock Computational Optimization and Applications 57(3):555–597. \url{https://doi.org/10.1007/s10589-013-9610-3}
	
	\bibitem[Aouad and Segev(2019)]{Aouad2019}
	Aouad~A, Segev~ D (2019)
	\newblock The ordered k-median problem: surrogate models and approximation
	algorithms.
	\newblock Mathematical Programming 177:55--83. \url{https://doi.org/10.1007/s10107-018-1259-3}
	
	
	\bibitem[Archetti et~al.(2016)Archetti, Corber\'{a}n, Plana, Sanchis, and
	Speranza]{Archetti16}
	Archetti~C, Corber\'{a}n~A, Plana~I, Sanchis~JM, Speranza~MG (2016)
	\newblock A branch-and-cut algorithm for the orienteering arc routing problem.
	\newblock Computers \& Operations Research 66:95--104. \url{https://doi.org/10.1016/j.cor.2015.08.003}
	
	
	\bibitem[Benati et~al.(2017)Benati, Puerto, and
	Rodríguez-Chía]{Benati2017}
	Benati~S, Puerto~J, Rodríguez-Chía~AM (2017)
	\newblock Clustering data that are graph connected.
	\newblock European Journal of Operational Research 261:43--53. \url{http://dx.doi.org/10.1016/j.ejor.2017.02.009}
	
	\bibitem[Blado and Toriello(2021)]{Blado2021}
	Blado~D, Toriello~A (2021)
	\newblock A column and constraint generation algorithm for the dynamic knapsack
	problem with stochastic item sizes.
	\newblock Mathematical Programming Computation 13:85--223. \url{https://doi.org/10.1007/s12532-020-00189-0}
	
	
	\bibitem[Boland et~al.(2006)Boland, Domínguez-Marín, Nickel, and
	Puerto]{Boland2006}
	Boland~N, Domínguez-Marín~P, Nickel~S, Puerto~J (2006)
	\newblock Exact procedures for solving the discrete ordered median problem.
	\newblock Computers \& Operations Research 33(11):3270--3300. \url{https://doi.org/10.1016/j.cor.2005.03.025}
	
	\bibitem[Calvino et~al.(2022)Calvino, López-Haro, Muñoz-Ocaña, Puerto, and
	Rodríguez-Chía]{Calvino2022}
	Calvino~JJ, López-Haro~M, Muñoz-Ocaña~JM, Puerto~J, 
	Rodríguez-Chía~AM (2022)
	\newblock Segmentation of scanning-transmission electron microscopy images
	using the ordered median problem.
	\newblock European Journal of Operational Research 302(2):671--687. \url{https://doi.org/10.1016/j.ejor.2022.01.022}
	
	\bibitem[CPLEX(2022)]{CPLEX}
	CPLEX (last accessed November 13, 2022)
	\newblock Differences between user cuts and lazy constraints.
	\newblock
	\href{https://www.ibm.com/docs/en/icos/20.1.0?topic=pools-differences-between-user-cuts-lazy-constraint}{Link to IBM page.} 
	
	\bibitem[de~Oliveira and Sagastizábal(2014)]{Oliveira2014}
	De~Oliveira~W, Sagastizábal~C (2014)
	\newblock Level bundle methods for oracles with on-demand accuracy.
	\newblock Optimization Methods and Software 29(6):1180--1209. \url{https://doi.org/10.1080/10556788.2013.871282}
	
	\bibitem[Deleplanque et~al.(2020)Deleplanque, Labbé, Ponce, and
	Puerto]{Deleplanque2020}
	Deleplanque~S, Labbé~M, Ponce~D, Puerto~J (2020)
	\newblock A branch-price-and-cut procedure for the discrete ordered median
	problem.
	\newblock INFORMS Journal on Computing 32:582--599. \url{https://doi.org/10.1287/ijoc.2019.0915}
	
	\bibitem[Deleplanque et~al.(2022)Deleplanque, Labbé, Ponce, and
	Puerto]{DOMPrepository}
	Deleplanque~S, Labbé~M, Ponce~D, Puerto~J (last accessed November 13, 2022)
	\newblock DOMP repository.
	\newblock 
	\href{https://gom.ulb.ac.be/gom/wp-content/uploads/2018/12/DOMP_Repository.zip}{Link to the repository.}
	
	\bibitem[Domínguez and Marín(2020)]{Dominguez2020}
	Domínguez~E, Marín~A (2020)
	\newblock Discrete ordered median problem with induced order.
	\newblock TOP 28:793--813. \url{https://doi.org/10.1007/s11750-020-00570-1}
	
	\bibitem[Edmonds and Johnson(1973)]{Edmonds}
	Edmonds~J, Johnson~EL (1973)
	\newblock Matching, Euler tours and the Chinese postman.
	\newblock Mathematical Programming 5(1):88--124. \url{https://doi.org/10.1007/BF01580113}
	
	
	\bibitem[Espejo et~al.(2009)Espejo, Marín, Puerto, and
	Rodríguez-Chía]{Espejo2009}
	Espejo~I, Marín~A, Puerto~J, Rodríguez-Chía~AM (2009)
	\newblock A comparison of formulations and solution methods for the
	minimum-envy location problem.
	\newblock Computers \& Operations Research 36(6):1966--1981. \url{https://doi.org/10.1016/j.cor.2008.06.013}
	
	\bibitem[Fernández et~al.(2017)Fernández, Pozo, Puerto, and
	Scozzari]{FernandezPozoPuertoScozzari}
	Fernández~E, Pozo~MA, Puerto~J, Scozzari~A (2017)
	\newblock Ordered weighted average optimization in multiobjective spanning tree
	problem.
	\newblock European Journal of Operational Research 260(3):886--903. \url{https://doi.org/10.1016/j.ejor.2016.10.016}
	
	
	\bibitem[Focacci et~al.(2002{\natexlab{a}})Focacci, Lodi, and
	Milano]{Focacci2002a}
	Focacci~F, Lodi~A, Milano~M (2002{\natexlab{a}})
	\newblock Embedding relaxations in global constraints for solving {TSP} and
	{TSPTW}.
	\newblock Annals of Mathematics and Artificial Intelligence 34:291--311 \url{https://doi.org/10.1023/A:1014492408220}
	
	\bibitem[Focacci et~al.(2002{\natexlab{b}})Focacci, Lodi, and
	Milano]{Focacci2002b}
	Focacci~F, Lodi~A, Milano~M (2002{\natexlab{b}})
	\newblock Optimization-oriented global constraints.
	\newblock Constraints 7:351--365. \url{https://doi.org/10.1023/A:1020589922418}
	
	\bibitem[Focacci et~al.(1999)Focacci, Lodi, Milano, and Vigo]{Focacci1999}
	Focacci~F, Lodi~A, Milano~M, Vigo~D (1999)
	\newblock Solving {TSP} through the integration of {OR} and {CP} techniques.
	\newblock Electronic Notes in Discrete Mathematics 1:3--25. \url{https://doi.org/10.1016/S1571-0653(04)00002-2}
	
	\bibitem[Fourour and Lebbah(2020)]{Fourour2020}
	Fourour~S, Lebbah~Y (2020)
	\newblock Equitable optimization for multicast communication.
	\newblock International Journal of Decision Support System Technology 12(3):1--25. \url{https://doi.org/10.4018/IJDSST.2020070101}
	
	\bibitem[Gleixner et~al.(2018)Gleixner, Bastubbe, Eifler, Gally, Gamrath,
	Gottwald, Hendel, Hojny, Koch, L{\"u}bbecke, Maher, Miltenberger, M{\"u}ller,
	Pfetsch, Puchert, Rehfeldt, Schl{\"o}sser, Schubert, Serrano, Shinano,
	Viernickel, Walter, Wegscheider, Witt, and Witzig]{GleixnerEtal2018ZR}
	Gleixner~A, Bastubbe~M, Eifler~L, Gally~T, Gamrath~G, Gottwald~RL,
	Hendel~G, Hojny~C, Koch~T, L{\"u}bbecke~ME, Maher~SJ,
	Miltenberger~M, M{\"u}ller~B, Pfetsch~ME, Puchert~C, Rehfeldt~D,
	Schl{\"o}sser~F, Schubert~C, Serrano~F, Shinano~Y, Viernickel~JM,
	Walter~M, Wegscheider~F, Witt~JT, Witzig~J (2018)
	\newblock {The SCIP Optimization Suite 6.0}.
	\newblock ZIB-Report 18-26, Zuse Institute Berlin.
	
	\bibitem[Gr{\"o}tschel and Holland(1985)]{Grtschel1985}
	Gr{\"o}tschel~M, Holland~O (1985)
	\newblock Solving matching problems with linear programming.
	\newblock Mathematical Programming 33:243--259. \url{https://doi.org/10.1007/BF01584376}
	
	\bibitem[Hakimi(1964)]{Hakimi1964}
	Hakimi~SL (1964)
	\newblock Optimum location of switching centers and the absolute centers and medians of a graph.
	\newblock Operations Research 12:450--459. \url{https://doi.org/10.1287/opre.12.3.450}
	
	\bibitem[Kalcsics et~al.(2010)Kalcsics, Nickel, Puerto, and Rodríguez-Chía]{Kalcsics2010}
	Kalcsics~J, Nickel~S, Puerto~J, Rodríguez-Chía~AM (2010)
	\newblock The ordered capacitated facility location problem.
	\newblock TOP 18:203--222. \url{https://doi.org/10.1007/s11750-009-0089-0}
	
	\bibitem[Labbé et~al.(2017)Labbé, Ponce, and Puerto]{Labbe2017}
	Labbé~M, Ponce~D, Puerto~J (2017)
	\newblock A comparative study of formulations and solution methods for the
	discrete ordered $p$-median problem.
	\newblock Computers \& Operations Research 78:230--242. \url{https://doi.org/10.1016/j.cor.2016.06.004}
	
	\bibitem[Letchford et~al.(2004)Letchford, Reinelt, and Theis]{Letchford}
	Letchford~AN, Reinelt~G, Theis~DO (2004)
	\newblock A faster exact separation algorithm for blossom inequalities.
	\newblock In: Bienstock~ D, Nemhauser~G (eds) Integer Programming
	and Combinatorial Optimization, Springer, Berlin Heidelberg, pp 196--205
	
	\bibitem[Magnanti and Wolsey(1995)]{Magnanti}
	Magnanti~TL, Wolsey~LA (1995)
	\newblock Chapter 9 optimal trees.
	\newblock In: Ball~MO, Magnanti~TL, Monma~CL, Nemhauser~GL (eds) Network Models, volume~7 of Handbooks in Operations Research and Management Science, Elsevier, pp 503--615
	
	\bibitem[Marín et~al.(2022)Marín, Martínez-Merino, Puerto, and
	Rodríguez-Chía]{Marin2022}
	Marín~A, Martínez-Merino~LI, Puerto~J, Rodríguez-Chía~AM (2022)
	\newblock The soft-margin support vector machine with ordered weighted average.
	\newblock Knowledge-Based Systems 237:107705. \url{https://doi.org/10.1016/j.knosys.2021.107705}
	
	\bibitem[Marín et~al.(2009)Marín, Nickel, Puerto, and Velten]{Marin2009}
	Marín~A, Nickel~S, Puerto~J, Velten~S (2009)
	\newblock A flexible model and efficient solution strategies for discrete
	location problems.
	\newblock Discrete Applied Mathematics 157(5):1128--1145. \url{https://doi.org/10.1016/j.dam.2008.03.013}
	
	\bibitem[Marín et~al.(2010)Marín, Nickel, and Velten]{Marin2010}
	Marín~A, Nickel~S, Velten~S (2010)
	\newblock An extended covering model for flexible discrete and equity location
	problems.
	\newblock Mathematical Methods of Operations Research 71(1):125--163. \url{https://doi.org/10.1007/s00186-009-0288-3}
	
	\bibitem[Marín and Pelegrín(2019)Marín, Pelegrín]{Marin2019}
	Marín~A, Pelegrín~M (2019)
	\newblock $p$-Median Problems
	\newblock In: Laporte, G., Nickel, S., Saldanha da Gama, F. (eds) Location Science. Springer, Cham, pp 25--50 \url{https://doi.org/10.1007/978-3-030-32177-2_2}
	
	\bibitem[Mart\'inez-Merino et~al.(2017)Mart\'inez-Merino, Albareda-Sambola, and
	Rodríguez-Chía]{MM2017}
	Mart\'inez-Merino~LI, Albareda-Sambola~M, Rodríguez-Chía~AM (2017)
	\newblock The probabilistic $p$-center problem: Planning service for potential customers.
	\newblock European Journal of Operational Research 262(2):509--520. \url{https://doi.org/10.1016/j.ejor.2017.03.043}	
	
	\bibitem[Mazzi et~al.(2021)Mazzi, Grothey, McKinnon, and Sugishita]{Mazzi2020}
	Mazzi~N, Grothey~A, McKinnon~K, Sugishita~N (2021)
	\newblock Benders decomposition with adaptive oracles for large scale
	optimization.
	\newblock Mathematical Programming Computation 13:683--703. \url{https://doi.org/10.1007/s12532-020-00197-0}
	
	
	\bibitem[Nickel(2001)]{Nickel2001}
	Nickel~S (2001)
	\newblock Discrete ordered Weber problems.
	\newblock In: Fleischmann~B, Lasch~R, Derigs~U, Domschke~W, Rieder~U (eds) Operations Research Proceedings, Springer, Berlin Heidelberg, pp 71--76
	
	\bibitem[Nickel and Puerto(2005)]{NickelPuerto2005}
	Nickel S, Puerto J (2005)
	\newblock Location Theory: a unified approach.
	\newblock Springer-Verlag, Berlin Heidelberg
	
	\bibitem[Ogryczak et~al.(2011)Ogryczak, Perny, and Weng]{Ogryczak2011}
	Ogryczak~W, Perny~P, Weng~P (2011)
	\newblock On minimizing ordered weighted regrets in multiobjective Markov
	decision processes, volume 6992 LNAI of Lecture Notes in Computer
	Science (including subseries Lecture Notes in Artificial Intelligence and
	Lecture Notes in Bioinformatics)
	
	\bibitem[Ponce et~al.(2018)Ponce, Puerto, Ricca, and Scozzari]{Ponce2018}
	Ponce~D, Puerto~J, Ricca~F, Scozzari~A (2018)
	\newblock Mathematical programming formulations for the efficient solution of
	the $k$-sum approval voting problem.
	\newblock Computers \& Operations Research 98:127--136. \url{https://doi.org/10.1016/j.cor.2018.05.014}
	
	\bibitem[Puerto et~al.(2013)Puerto, Ramos, and Rodríguez-Chía]{Puerto2013}
	Puerto~J, Ramos~AB, Rodríguez-Chía~AM (2013)
	\newblock A specialized branch \& bound \& cut for single-allocation ordered
	median hub location problems.
	\newblock Discrete Applied Mathematics 16:2624--2646. \url{https://doi.org/10.1016/j.dam.2013.05.035}
	
	\bibitem[ReVelle and Swain(1970)]{ReVelle1970}
	ReVelle~CS, Swain~R (1970)
	\newblock Central facilities location.
	\newblock Geographical Analysis 2:30--42. \url{https://doi.org/10.1111/j.1538-4632.1970.tb00142.x}
	
	\bibitem[SCIP(2022a)]{SCIP1}
	SCIP (last accessed November 13, 2022a)
	\newblock When should I implement a constraint handler, when should I implement a separator?
	\newblock \url{https://www.scipopt.org/doc/html/FAQ.php#conshdlrvsseparator}
	
	\bibitem[SCIP(2022b)]{SCIP2}
	SCIP (last accessed November 13, 2022b)
	\newblock SCIPcreateCons().
	\newblock
	\href{https://www.scipopt.org/doc/html/group__PublicConstraintMethods.php#ga8d771b778ea2599827e84d01c0aceaf2}{Link to SCIP page.}

	
	\bibitem[Tamir(2001)]{Tamir2001}
	Tamir.~A (2001)
	\newblock The $k$-centrum multi-facility location problem.
	\newblock Discrete Applied Mathematics 109(3):93--307. \url{https://doi.org/10.1016/S0166-218X(00)00253-5}
	
	\bibitem[Wolf et~al.(2014)Wolf, Fábián, Koberstein, and Suhl]{Wolf2014}
	Wolf~C, Fábián~CI, Koberstein~A, Suhl~L (2014)
	\newblock Applying oracles of on-demand accuracy in two-stage stochastic
	programming – A computational study.
	\newblock European Journal of Operational Research 239(2):437--448. \url{https://doi.org/10.1016/j.ejor.2014.05.010}
	
\end{thebibliography}
\end{document}